\documentclass[11pt,a4paper]{article}
\usepackage[utf8]{inputenc}
\usepackage[english]{babel}
\usepackage[left=1.5cm,right=1.5cm,top=2cm,bottom=2cm,foot=1cm]{geometry}
\usepackage{amsmath,amssymb,latexsym, amsthm}
\renewcommand\footnotemark{}
\newtheorem{thm}{Theorem}

\newtheorem{prop}{Proposition}

\theoremstyle{definition}

\title{Diffeomorphism Groups of Compact 4-manifolds are not always Jordan }
\begin{document}
\author{Bal\'azs Csik\'os,
L\'aszl\'o Pyber, Endre Szab\'o
} \date{}
\maketitle
\begin{abstract} We show that if $M$ is a compact smooth manifold  diffeomorphic to the total space of an orientable $S^2$ bundle over the torus $T^2$, then its diffeomorphism group does not have the Jordan property, i.e., Diff$(M)$ contains a finite subgroup $G_n$ for any natural number $n$ such that every abelian subgroup of $G_n$ has index at leat $n$. This gives a counterexample to an old conjecture of Ghys.
\end{abstract}

\maketitle

\section{Introduction} By Jordan's classical theorem (see, e.g., \cite{Wolf}), for any natural number $n$, there is a constant $c_n$ depending only on  $n$, such that every finite subgroup of the group $GL(n,\mathbb C)$ has an abelian subgroup of index at most $c_n$.

We say that a group $G$ has the \emph{Jordan property} if there is a constant $c$ depending only on the group $G$, such that every finite subgroup of $G$ has an abelian subgroup of index at most $c$.

Motivated by the Jordan theorem, E. Ghys conjectured that the diffeomorphism groups of compact smooth manifolds have the Jordan property. As noted in \cite{Riera}, this conjecture was discussed in several talks of Ghys \cite{Ghys}, but appeared in print for the first time in \cite{Fisher}. 

Known examples of non-compact manifolds \cite{Popov}, \cite{Zimmermann} having diffeomorphism groups  which do not possess the Jordan property show that the compactness assumption in the conjecture is essential. 

The conjecture of Ghys was verified in several special cases. Zimmermann \cite{Zimmermann} proved the conjecture for compact $3$-manifolds and in a series of papers \cite{Riera0}, \cite{Riera},  \cite{Riera2}, \cite{Riera3}, Riera proved the conjecture for tori, acyclic manifolds, homology spheres, and manifolds with non-zero Euler characteristic.

The goal of this paper is to prove the following theorem.
\begin{thm} $M$ is diffeomorphic to $T^2\times S^2$ or to the total space of a nontrivial smooth orientable $S^2$ bundle over $T^2$, then the diffeomorphism group of $M$ does not have the Jordan property.
\end{thm}
The proof uses some standard facts on smooth complex line bundles and sphere bundles over the torus, we summarize in section 2. 

Then we apply ideas from algebraic geometry that appeared in Yu.G. Zarhin \cite{Zarhin}, where the same ideas were used to  prove that the group of birational automorphisms of the product of an elliptic curve and a projective line over an algebraically closed field of characteristic zero does not have the Jordan property. We collect the necessary information on holomorphic line bundles and prove the main theorem in section 3.

\section{Smooth complex line bundles and sphere bundles over the torus} 

Let $T^2=\mathbb R^2/\mathbb Z^2$ be a $2$-dimensional torus with a given complex structure, and fix a point $o\in T^2$. As the universal covering space of $T^2$ is biholomorphically equivalent to $\mathbb C$, we can identify $T^2$ with the factor space $\mathbb C/\Gamma$, where $\Gamma$ is a rank $2$ lattice in $\mathbb C$ acting on $\mathbb C$ by translations. Addition in $\mathbb C$ induces a commutative group operation $+$ on $T^2$ with neutral element $o$. For $x\in T^2$, we shall denote by $T_x\colon T^2\to T^2$, $y\mapsto x+y$ the translation by $x$.

Smooth complex line bundles $\xi=(E\xrightarrow{\pi}T^2)$ over $T^2$ are classified by the first Chern class $c_1(\xi)\in H^2(T^2,\mathbb Z)\cong \mathbb Z$, or, equivalently, by the Chern number $\int_{T^2}c_1(\xi)\in \mathbb Z$, where the integral is computed using the orientation of $T^2$ induced by the complex structure. Denote by $\xi_n=(E_n\xrightarrow{\pi_n}T^2)$ the complex line bundle with Chern number $n\in \mathbb Z$. It is known that $\xi_n\otimes_{\mathbb C}\xi_m\cong \xi_{n+m}$, in particular $\xi_n$ is the $n$th tensor power of $\xi_1$ for all $n>0$, $\xi_0$ is the trivial bundle.

We can associate to any complex line bundle $\xi$ a smooth complex projective line bundle $P(\xi\oplus \xi_0)$ with fibers diffeomorphic to $\mathbb CP^1\cong S^2$. Denote by $Y_n$ the total space of the bundle $P(\xi_n\oplus \xi_0)$. It is known that non-isomorphic sphere bundles over a surface can have diffeomorphic total spaces. According to Theorems 1 and 2 in \cite{Melvin}, the total spaces of two $S^2$-bundles $\eta_1$ and $\eta_2$ over $T^2$ are diffeomorphic, if and only if the  Stiefel-Whitney classes $w_1$ and $w_2$ of them satisfy the following conditions:
\begin{itemize}
	\item[(i)] $w_2(\eta_1)=w_2(\eta_2)$;
	\item[(ii)] either $w_1(\eta_1)=w_1(\eta_2)=0$ or none of $w_1(\eta_1)$ and $w_1(\eta_2)$ is $0$.
\end{itemize}
As a corollary, the total spaces of all smooth sphere bundles over $T^2$ belong to one of four different diffeomorphism classes. Due to the presence of complex structures, the sphere bundles $P(\xi_n\oplus \xi_0)$ are all orientable, so their first Stiefel-Whitney classes vanish. The second Stiefel-Whitney class $w_2(P(\xi_n\oplus \xi_0))\in H^2(T^2,\mathbb Z_2)\cong \mathbb Z_2$ is the mod 2 reduction of $n$, so we have the following
\begin{prop}
$Y_n$ is diffeomorphic to $Y_m$ if and only if $n\equiv m$ (mod $2$).
\end{prop} 
\section{Holomorphic structures on complex line bundles and their biholomorphic automorphism groups}
Consider the complex line bundle $\xi_1$. There are infinitely many holomorphic structures on $\xi_1$ compatible with the given complex structure on $T^2$. Let us fix any of these holomorphic structures and equip $\xi_n$ with the holomorphic structure for which the isomorphisms $\xi_n\otimes_{\mathbb C}\xi_m\cong \xi_{m+n}$ become holomorphic isomorphism. 

Mumford \cite{Mumford} considered the subgroup  
$$
H(\xi_n)=\{x\in T^2 : T_x^*(\xi_n)\cong \xi_n\}
$$
of $T^2$, where $\cong$ means the isomorphism of \emph{holomorphic} line bundles. It is known that $H(\xi_n)$ is finite if and only if $\xi_n$ is an ample line bundle, which is the case if and only if  $n>0$. 

There is an extension of the group $H(\xi_n)$ defined by
$$
\mathcal G(\xi_n)=\{(x,\phi):x\in H(\xi_n), \text{ $\phi\colon E_n\to E_n$ is a biholomorphic map such that $\pi_n\circ \phi=T_x\circ \phi$}\}.
$$
The kernel of the natural homomorphism $\mathcal G(\xi_n)\to H(\xi_n)$, $(x,\phi)\mapsto x$ is $\mathbb C^*$ so we have a short exact sequence
$$
0\to \mathbb C^*\to \mathcal G(\xi_n)\to H(\xi_n)\to 0.
$$

The key observation is that every holomorphic structure on the bundle $\xi_n$ induces a unique holomorphic structure on the associated projective line bundle $P(\xi_n\oplus \xi_0)$, hence a complex structure on $Y_n$. For any $(x,\phi)\in \mathcal G(\xi_n)$, $\phi$ extends uniquely to a biholomorphic map $\tilde \phi\colon Y_n\to Y_n$. The map $(x,\phi)\mapsto \tilde \phi$ provides an embedding of the group $\mathcal G(\xi_n)$ into the diffeomorphism group of $Y_n$. Combining this fact with Proposition 1, we obtain the following 
\begin{prop}
The diffeomorphism group of $Y_n$ contains subgroups isomorphic to $\mathcal G(\xi_m)$ for every $m>0$ satisfying $m\equiv n$ mod $2$.
\end{prop}
The structure of the groups $ H(\xi_n)$ and $\mathcal G(\xi_n)$ was described in details by Mumford \cite{Mumford}. Recall those elements of the description that are relevant for the proof of our main result. 
\begin{thm}[Mumford, \S 1 in \cite{Mumford}]
The group $H(\xi_n)$ contains an (abelian) subgroup $K$ such that $H(\xi_n)$ is isomorphic to the group $K\oplus \hat K$, where $\hat K$ is the multiplicative group of characters $K\to S^1\subset \mathbb C^*$.

The group $\mathcal G(\xi_n)$ is isomorphic to the group defined on the set $\mathbb C^*\times (K\oplus \hat K)$ with the multiplication rule
$$
(a,k\oplus l)\cdot(a',k'\oplus l')=(a\cdot a'\cdot l'(k), (k+k')\oplus(l\cdot l')).
$$
\end{thm}

\begin{thm}[Mumford, Prop. 4 in \S 2 of \cite{Mumford}]
If $k$ is a positive integer, then  
$$
H(\xi_{kn})=\{x\in T^2 : kx\in H(\xi_{n})\}.
$$
In particular, the case $n=1$ yields that $H(\xi_{k})$ has at least $k^2$ elements.
\end{thm}

We obtain as a corollary that if $H(\xi_n)$ has $N$ elements, then $N\geq n^2$, $N$ is a square number, $K$ has $\sqrt{N}$ elements. Moreover, if $\mathbb Z_{\sqrt{N}} <\mathbb C^*$ denotes the cyclic subgroup of  order $\sqrt{N}$, then the subset $G_n=\mathbb Z_{\sqrt{N}}\times (K\oplus \hat K)\subset\mathbb C^*\times (K\oplus \hat K)$ is a finite subgroup of $\mathcal G(\xi_n)$.

Now we are ready to prove our main theorem.

\begin{proof}[Proof of Theorem 1.]
As it was pointed out Y.G. Zarhin \cite{Zarhin}, the index of any abelian  subgroup of $G_n$ is at least $\sqrt{N}\geq n$. By Proposition 1, $M$ is diffeomorphic to $Y_m$, where $m$ is either $0$ or $1$. In both cases, Diff$(M)$ contains an infinite sequence of finite subgroups $G_n$, ($n>0$ and $n\equiv m$ mod $2$), having the property that any abelian subgroup of $G_n$ has index at least $n$. This proves the theorem.  
\end{proof}
\section{Acknowledgments}

The first author is supported by the Hungarian National Science and Research Foundation OTKA K112703. During the research he also enjoyed the hospitality of the Alfr\'ed R\'enyi Institute of Mathematics as a guest researcher. 

The second author is supported in part by K84233.

The third author is supported in part by OTKA NK81203, K84233,
and   by MTA Rényi "Lendület" Groups and Graphs Research Group.

The authors are indebted to Andr\'as N\'emethi and Andr\'as Stipsicz for helpful and fruitful discussions.
\bibliographystyle{ieeetr}
\bibliography{Ghys}
Bal\'azs Csik\'os, Institute of Mathematics, E\"otv\"os Lor\'and University, Budapest, P\'azm\'any P. stny. 1/C, H-1117 Hungary. \emph{E-mail address:} csikos@cs.elte.hu

L\'aszl\'o Pyber, Alfr\'ed R\'enyi Institute of Mathematics, Budapest, Re\'altanoda u. 13-15, H-1053 Hungary, \emph{E-mail address:} pyber.laszlo@renyi.mta.hu

Endre Szab\'o, Alfr\'ed R\'enyi Institute of Mathematics, Budapest, Re\'altanoda u. 13-15, H-1053 Hungary, \emph{E-mail address:} szabo.endre@renyi.mta.hu
\end{document}